\documentstyle[12pt]{article}
\begin{document}

\newtheorem{theorem}{Theorem}[section]  
\newtheorem{definition}[theorem]{Definition}  
\newtheorem{example}[theorem]{Example}  
\newtheorem{examples}[theorem]{Examples}
\newtheorem{lemma}[theorem]{Lemma}  
\newtheorem{proposition}[theorem]{Proposition}  
\newtheorem{corollary}[theorem]{Corollary}  
\newtheorem{remark}[theorem]{Remark}  
\newtheorem{conjecture}[theorem]{Conjecture}

\newcommand{\A}{{\cal A}} 
\newcommand{\B}{{\cal B}} 
\newcommand{\C}{{\cal C}}
\newcommand{\D}{{\cal D}}
\newcommand{\E}{{\cal E}}
\newcommand{\F}{{\cal F}}
\newcommand{\Ho}{{\cal H}}
\newcommand{\Io}{{\cal I}}
\newcommand{\K}{{\cal K}}
\newcommand{\Lo}{{\cal L}}
\newcommand{\M}{{\cal M}}
\newcommand{\T}{{\cal T}}
\newcommand{\SS}{{\cal S}}
\newcommand{\X}{{\cal X}}
\newcommand{\al}{\alpha} 
\newcommand{\Ann}{\mbox{\rm Ann}}

\newcommand{\cn}{{\bf {\rm C}} 
\hspace{-.4em}      {\vrule height1.5ex width.08em depth-.04ex} 
\hspace{.3em}} 
\newcommand{\cp}{\cn{{\bf {\rm P}}
\hspace{-.4em}      {\vrule height1.5ex width.08em depth-.04ex} 
\hspace{.3em}}}
\newcommand{\pr}{{{\bf {\rm P}}
\hspace{-.4em}      {\vrule height1.5ex width.08em depth-.04ex} 
\hspace{.3em}}}
\newcommand{\codim}{\mbox{\rm codim}}
\newcommand{\depth}{\mbox{\rm depth}}
\newcommand{\Der}{\mbox{\rm Der}(S)} 
\newcommand{\gl}{{g_{\lambda}}} 
\newcommand{\Gr}{\mbox{\rm Gr}} 
\newcommand{\Hom}{\mbox{\rm Hom}}
\newcommand{\im}{\mbox{\rm Im}}
\newcommand{\ints}{{\sf Z}\hspace{-.36em}{\sf Z}} 
\newcommand{\complex}{{{\sf C} \kern -.53em {\rm l} \kern +0.38em}}
\newcommand{\rats}{{{\sf Q} \kern -.45em {\rm l} \kern +.30em}} 
\newcommand{\reals}{{{\rm l} \kern -.15em {\sf R} }}	
\newcommand{\la}{\lambda}
\newcommand{\om}{\omega}  
\newcommand{\Om}{\Omega}  
\newcommand{\Op}{\Omega^{p}}  
\newcommand{\OpA}{\Omega^{p}(\A)}

\newcommand{\p}{\partial}  
\newcommand{\PD}{\mbox{\rm pd}}
\newcommand{\Pl}{{\Phi_{\lambda}}}  
\newcommand{\pl}{{\phi_{\lambda}}}  
\newcommand{\Poin}{\mbox{\rm Poin}} 
\newcommand{\PoinO}{\Poin(\Omega^{*}(c\A); u, } 
\newcommand{\PoinGr}{\Poin(\Gr\Omega^{*}(\A); u, } 
\newcommand{\rk}{\mbox{\rm rk}}
\newcommand{\cd}{\mbox{\rm cd}}
\newcommand{\sign}{\mbox{\rm sign}}
\newcommand{\proof}{{\bf Proof.~}} 
\newcommand{\qed}{~~\mbox{$\Box$}} 

\newcommand{\ra}{{\rightarrow}}  
\newcommand{\rn}{{\rm I}\hspace{-.2em}{\rm R}}

\newcommand{\scn}{\scriptsize\cn} 
\newcommand{\stR}[1]{\stackrel{#1}{\longrightarrow}} 
\newcommand{\th}{\theta}

\newcommand{\we}{\wedge}
\newcommand{\ar}{\buildrel D_{\la}\over\rightarrow}

\title{Small rational model of subspace complement}
\author { 
{\sc Sergey Yuzvinsky }\\
{\small\it University of Oregon,
Eugene, OR 97403 USA}\\
{\small\it yuz@math.uoregon.edu} }

\maketitle

\section{Introduction}
\bigskip
This paper concerns the cohomology ring of the complement of a complex subspace
arrangement. Let $\A$ be a finite set of proper subspaces of a finite
dimensional complex linear space $V$.
Put $C(\A)=V\setminus\bigcup_{A\in \A}A$. We can assume that there are no inclusions
among elements of $\A$. A rough characterization
 of $\A$ is given by the lattice $X=X(\A)$
consisting of all the intersections $A$
of the elements from $\A$ ordered opposite to
inclusion and labeled by their codimensions $\cd A$. 

 If each subspace $A\in\A$ is a hyperplane then $\A$ can be given by a set
$\{\alpha_A\in V^*\}$ where $\ker\alpha_A=A$. In this case, the lattice $X$ is
geometric and the codimension of its each element can be recovered from the ordering
as the rank of this element.
There are two main results describing the structure of the ring
$R=H^*(C(\A),{\ints})$ in
this case.
The first one is the Arnold-Brieskorn theorem \cite{Ar,Br} stating that the algebra of
differential forms generated by the closed forms ${1\over{2\pi i}}{{d\alpha_A}\over
{\alpha_A}}$ 
is isomorphic to $R$ under the de Rham homomorphism.
The second one is the Orlik-Solomon theorem proving  that $R$ depends only on
$X$. Moreover the theorem gives a presentation of $R$ by generators and
relations
defined by $X$.

In the general case, $X$ is not necessarily geometric
and the codimensions ${\rm cd}A$ of $A\in X$ 
can not be anymore recovered from the order on
$X$. Again, there are two main results about $H^*(C(\A))$.
The first one is the Goresky-MacPherson theorem \cite{GMP} proving that the groups
$H^p(C(\A))$ are defined by the labeled lattice
 $X$. More precisely these groups are expressed 
in terms of the local homology groups
 of $X(\A)$ as follows.
For every $p$
$$\tilde H^p(C(\A),\ints)=\oplus_{A\in X}\tilde H_{2{\rm cd}A-p-2}((V,A),\ints)$$
where $(V,A)=\{B\in X|V<B<A\}$. 
 This is an example of application of
Stratified Morse theory developed in \cite{GMP}.
 Later, this result was obtained by 
a different method in \cite{ZZ}. A completely different approach was used recently
by De Concini and Procesi \cite{DCP}.
They constructed
a rational model for $C(\A)$ using only the labeled lattice $X$
 (see also section 2) and proved that
the rational cohomology algebra and rational 
homotopy type of $C(\A)$ are defined by this lattice.
The natural problem that was left open is how to recover the ring structure explicitly
from the combinatorics, in particular how to relate the multiplication
to the local homology of $X$ from the Goresky-MacPherson
formulas.

The present paper solves this problem (Theorem 6.5).
To be more specific here we need to introduce more notation. 
For every pair $(A,B)$ ($A,B\in X$)
 denote by $A\vee B$ the join (i.e., the least upper bound) 
of the pair. The joins of sets of atoms from $(V,A)\cup (V,B)$ 
form a subposet $X_{A,B}$ of $(V,A\vee B)$. The following condition on $A,B$ is 
an important breaking point
$${\rm cd}(A\vee B)={\rm cd} A+{\rm cd} B.\eqno(*)$$
If this condition holds then the flag complex of $X_{A,B}$ has the homotopy type of
the suspension of the join of the flag complexes of $(V,A)$ and $(V,B)$.
Then there is a natural isomorphism
$$\phi_{A,B}:\tilde H_r((0,A),{\rats})\otimes \tilde H_s((0,B),{\rats})\to
\tilde H_{r+s+2}(X_{A,B},{\rats})$$
(e.g., see \cite{Mu}).

{\bf Theorem}.
{\it The multiplication in $H^*(C(\A),\rats)$ is given by  bilinear pairings

$$\psi_{A,B}:\tilde H_r((0,A),{\rats})\otimes \tilde H_s((0,B),{\rats})\to
\tilde H_{r+s+2}((0,A\vee B),{\rats})$$
for every $A,B\in X$. Here $\psi_{A,B}$
is the composition of $\phi_{A,B}$ with the
embedding $X_{A,B}\subset (V,A\vee B)$ if (*) holds and $\psi_{A,B}=0$ otherwise.
}

Notice that the set of flag complexes of $X(\A)$ for all the subspace arrangements
$\A$ includes the homotopy types of all finite simplicial
complexes. In this generosity, it is impossible to give any explicit
description of local homology classes of $X(\A)$ and the multiplication on
them. Instead we describe the multiplication on the local flag complexes
that induces the right multiplication on the homology.

 Let $A,B\in X(\A)$ with condition (*). Let
 $F_A=(A_1<A_2<\cdots<A_i)$ and $F_B=(B_1<B_2<\cdots <B_j)$ be
flags in the open intervals
$(V,A)$ and $(V,B)$ respectively. Augment the flags by $A$ and $B$
respectively and take the shuffle product of them.
The result is a linear combination (with coefficients equal to signs of
permutations) of certain sequences of $A_r$ and $B_s$. Create a flag from each
such sequence taking the intersections of all initial subsequences. 
 Equate flags with repetitions to 0 and delete $A\vee B$ at the ends
of all others. The resulting linear combination of flags is the product
$F_A\cdot F_B$.
This product generates the product on linear combinations of flags 
 by bilinearity that commutes with the differentials. Hence
it induces a bilinear pairing
$$\tilde H_{i-1}((V,A))\otimes \tilde H_{j-1}((V,B))\to 
\tilde H_{i+j}((V,A\vee B)).$$
This is the multiplication on local homology of $X$.
The multiplication can be defined even easier on so called Whitney (or relative
atomic) complexes but we will present it elsewhere.

 For the classes of lattices whose
local homology can be described explicitly our theorem gives presentations
of $H^*(C(\A)$. We consider two  such classes; first one where
$X$ is
a geometric lattice and the other one where $C(\A)$ is a so called
$k$-equal manifold. In the former case, the local homology groups of $X$ 
were computed in
\cite{Fo,OT}; in the latter case, they were computed 
in \cite{BW}. 
We give presentations of the algebra $H^*(C(X),\rats)$ in
these cases that  in particular solves 
a problem stated in \cite{BW}.
The representation for the former case was obtained independently and by completely
different method by E.M. Feichtner in \cite{Fe}.

The method we use  in this paper is as follows.
We start with the De Concini - Procesi differential graded
algebra $M$ that is a specialization of the Morgan rational model
of the complement of a divisor with normal crossings (see \cite{Mo}).
The algebra $M$ is 
a rational model of $C(\A)$. We find inside $M$
a significantly smaller subalgebra $CM$ quasi-isomorphic to $M$ whence 
also a rational model of $C(\A)$ (that may be of interest by itself).
The algebra $CM$ gives a multiplicative structure on the flag complexes of $X$
that
induces the ring structure on  $H^*(C(X),\rats)$ described above.
Moreover $CM$ has a natural integral version $CM(\ints)$
that gives the right groups $H^p(C(X))$
 and is a good candidate for an integral model of
$C(X)$ (see Remark \ref {ints} and
Conjecture \ref{int}). 

The paper is organized as follows. In section 2, we recall the 
De Concini-Procesi model $M$ 
and a monomial basis of it found in \cite{Yu}.
 In section 3, we switch to generators more convenient for
the purpose of the paper and exhibit a basis of monomials
in the new generators.
 In section 4, we define the smaller subcomplex $CM$ that is
a deformation retract of $M$.
In section 5, we prove that $CM$ is a subalgebra of $M$ and
study the multiplicative structure on it.
In section 6, we describe the multiplication induced on its cohomology.
In section 7 and 8, we give a presentation
of the algebra
$H^*(C(X))$
for the cases where $X$ is a geometric lattice and where
$C(X)$ is a  $k$-equal manifold.

This work was started when I was visiting Technische Universit\"at in Berlin. I
want to express my gratitude to G. Ziegler for his hospitality. I am also
grateful to C. De Concini for useful communications and to Eva Feichtner for finding
an error in the first version of Proposition 7.2..
\bigskip
\section{De Concini - Procesi model}
\bigskip
In this section, we recall the construction of rational models
from \cite{DCP}.

To follow closer the notation in \cite{DCP} we dualize our initial definitions from
Introduction.
 Let $V$ be a finite dimensional complex linear space and $\A$ a finite set of
nonzero linear
subspaces of $V$. Without any loss of generality we can assume that their 
sum is $V$. The topological space we are interested in is 
$C(\A)=V^*\setminus\bigcup_{A\in
\A}A^{\circ}$ where $A^{\circ}$ is the annihilator of $A$ in $V^*$. The space
$C(\A)$ does not change 
if we include all nonempty sums of elements of $\A$ in $\A$.
Therefore we will be dealing with the lattice $X=X(\A)$ of all sums of elements of
$\A$ ordered by inclusion and labeled by the dimensions of its elements. 
Notice that $C(X)=C(\A)=V^*\setminus \bigcup_{A\in X,A\not=0}A^{\circ}$.
Let us emphasize that the only structure on $X$ used in this
paper is that of a labeled lattice.
We also denote by ${\bf P}(A)$ the projectivization of
a linear space $A$ and put ${\bf P}C(X)={\bf P}V^*\setminus\bigcup_{A\in X,A\not=0}
{\bf P}(A^{\circ})$. 

Now we describe two differential graded algebras ( DGA) $M^*(X)$ and $M(X)$.
 There are different (though quasi-isomorphic to each other)
versions of these algebras
depending on the choice of a building set of $X$.
 Recall that a subset $G$ of $X$
is a building set if every  $A\in X$ is the 
direct sum of maximal elements of $G$ contained in $A$. For the most part
of the paper we take the whole $X$ as 
the building set. The only time a different building set is
used is in the paragraph following Theorem \ref{DCP}.

 Let $G$ be a building set containing $V$.
Start with the polynomial algebra $P$ over ${\rats}$ with indeterminates
$c_A$ corresponding to $A\in G$
 and the exterior algebra $\Lambda$ over $\rats$ generated by
$e_A$ corresponding to 
$A\in G\setminus \{V\}$. Consider the DGA $P\otimes \Lambda$ with the grading
induced by $\deg c_A=2$ and $\deg e_A=1$ and the differential generated by
$d(e_A)=c_A$. Now consider the elements 
$$r(X_1,X_2,B)=
\prod_{A\in X_1}e_A\prod_{A\in X_2}c_A(\sum_{C\supseteq B}c_C)^{d(Y,B)}\eqno
(2.1)$$
of $P\otimes\Lambda$ 
where $X_1,X_2\subseteq G\setminus\{V\}$, $Y=X_1\cup X_2$, 
$B\in G,  \ B\supset A$ for all $A\in Y$ and $d(Y,B)=\dim B-\dim(\sum_{A\in Y}
A)$. Notice that in order to define the first product and thus $r(X_1,X_2,B)$
precisely one needs to introduce a linear order on $X_1$. Otherwise these
elements are defined only up to $\pm 1$ that suffices for the definition.
Let $J$ be the ideal generated by all the elements $r(X_1,X_2,B)$.
Clearly it is homogeneous and invariant with respect to $d$. Thus $M^*(X)
=(P\otimes \Lambda)/J$ is a DGA. DGA $M(X)$ is defined similarly with the
only difference that both $e_A$ and $c_A$ are indexed by all the elements
of $G$ (including $V$).
We will denote by $M^*(X)_r$ and $M(X)_r$
 the homogeneous components of degree $r$ of the algebras.

As it is customary,
 we will keep the notation $c_A$ and $e_A$ for the images of
these elements under factorization and the term ``monomials''
 for the images of monomials in $c_A$ and $e_A$.

\begin{theorem}(\cite{DCP}, Theorem 5.3.)
\label{DCP}
There is a $\rats$-algebra isomorphism
$\phi: H^*(M^*(X))\approx H^*({\bf P}C(X),{\rats})$. 
\end{theorem}

\medskip
In \cite{DCP}, the problem of the naturality of $\phi$ is left
open. It follows from \cite{Mo} that after tensoring with $\cn$ this isomorphism
becomes canonical. Over $\rats$ it is not functorial in the category of 
smooth algebraic varieties. It is probably functorial in the category of
arrangement complements but a proof is not known to me.

Now we want a model for $C(X)\subset V^*$.
 We can view $X$
 as a
subspace arrangement in $V'=V\oplus W$ with $\dim W=1$ (this is a
generalization of the coning construction from
 \cite{OT}, p.14).
Then we can
consider the new arrangement $X'=X\cup\{W\}\cup\{A\oplus W|A\in X\}$ in
$V'$. Notice that
 ${\bf P}C(X')=C(X)$.
 Thus $M^*(X')$ can serve as a model of $C(X)$.
To simplify definition of $M^*(X')$ 
we choose $X\cup\{W,V'\}$
as the building set for $X'$. Then to obtain $M^*(X')$ one needs to adjoin 
$e_W,e_V,c_{V'}$ and $c_W$ to the generators of $M^*(X)$. Notice that 
since $\dim W=1$ we have the relation
$$c_W+c_{V'}=0.$$
Thus $c_W=-c_{V'}$ and we can forget about $c_W$.
If we put
 $e_{V'}=-e_W$  then we see that $M^*(X')=M(X'')$ where
 $X''=X\cup \{V'\}$ is a subspace arrangement in $V'$.
The final comment is that $M(X)$ is
quasi-isomorphic to $M(X'')$ and thus can be taken as a model of $C(X)$.
This can be proved either by exhibiting an explicit homotopy equivalence
or by noticing
that the local homology of $X''$ does not give anything new with respect to
that of $X$ (cf. Section 6). 
Keeping this in mind we will focus our attention on $M=M(X)$ from now on.
Let us emphasize again that the building set for $M$ is the whole $X$.

Our goal for the rest of the section is to exhibit a linear basis of $M$
consisting of monomials in the generators $c_A$ and $e_A$. For each
$T\subseteq X$ denote by $\mu(T)$ the monomial of $\Lambda$ corresponding
to $T$ (defined up to $\pm$). Then consider
the grading $\Lambda=\oplus_{T\subseteq X}\Lambda^T$ of $\Lambda$
where $\Lambda^T={\bf C} \mu(T)$.
It generates a grading 
on $\Lambda\otimes P$ and, since every $r(X_1,X_2,B)$ is
homogeneous in this grading, a grading $M=\oplus_{T\subseteq X}M^T$ on $M$.
Then notice
that for every $A_1,A_2\in X$ neither of which is a subset of the other
 $$a_{A_1}a_{A_2}=0$$
in $M$ where $a_A$ denotes either $c_A$ or $e_A$. To see this it is enough
to take in (2.1) $Y=\{A_1,A_2\}$ and $B=A_1+ A_2$ that gives
$d(Y,B)=0$. Hence all the monomials in $e_A$ identified up to the ordering (or
 $\pm 1$)
are in one-to-one correspondence
with flags (by inclusion) in $X$. In particular, in the representation
$M=\oplus_{T\subseteq X}M^T$, a summand $M^T$ is nonzero only if $T$ is a flag
of $X$. In any case $M^T$ has the following
 structure of a graded algebra
(although not a subalgebra of $M$). 
Denote by $X^T$ the set of all elements of $X$ each one of which forms a
flag with $T$. Then $M^T$ 
is the factor of the polynomial ring ${\rats}[c_A|A\in X^T]$ by the
ideal generated by $$r(Y,B)=\prod_{A\in Y}c_A(\sum_{C\supseteq B}c_C)^{d^T(Y,B)}
$$
for every $Y\subset X^T\setminus T$, $B\in X^T$ such that $B\supset A$
 for every $A\in Y$, and $d^T(Y,B)=\dim
B-\dim(\sum_{A\in Y\cup T,A\subset B}A)$. (In fact in \cite{DCP} the algebras
$M^T$ appeared first and $M^*(X)$ 
was constructed from them following \cite{Mo}.)

A monomial basis of algebra $M^T$ was first constructed in \cite{Yu} (see
\cite{Yu}, Remark 3.11) and then generalized in \cite{Ga}. Combining the
monomial bases of $M^T$ for all flags $T$ we obtain the following.

\begin{proposition}
\label{basis}
Suppose $S,T\subseteq X$,
$T=\{A_1,\ldots,A_k\}$ is ordered by inclusion 
and $S\cup T$ is a flag. 
Suppose $m:S\to {\bf Z}_+$ 
satisfies $m(A)<\dim A-\dim A'$ for every $A\in S$ where $A'=0$ if $A$ is
the smallest element of $T\cup S$ and $A'$ is the predecessor of $A$ in $T\cup
S$ otherwise. 
Then the monomials $\lambda(S,T,m)=e_{A_1}\cdots e_{A_k}\prod_{B\in
S}c_B^{m(B)}$ form a basis of $M$. 
\end{proposition}

\medskip
We will need later
 one more idea from \cite{Yu}. Consider an arbitrary monomial
$\lambda=e_{A_1}\cdots e_{A_k}c_{B_1}^{m_1}\cdots c_{B_{\ell}}^{m_{\ell}}$.
Then define the {\it weight} of $\lambda$ as
$${\rm wt}(\lambda)=\sum_{i=1}^k\dim A_i+2\sum_{i=1}^{\ell}m_i\dim B_i.$$

It is easy to see that for any basic monomial $\lambda(S,T,m)$ that
appears in the
linear decomposition of $\lambda$ with a nonzero coefficient we have
$${\rm wt}(\lambda(S,T,m))\leq {\rm wt}(\lambda)$$
 (cf. \cite{Yu}, Proof of Proposition 2.2).
\bigskip
\section{New generators}
\bigskip
In this section, we switch to new generators of algebra $M$.
 For every $A\in X$ put
$$\sigma_A=\sum_{B\supseteq A}c_B\ {\rm and}\ \tau_A=\sum_{B\supseteq A}e_B.$$
 Clearly $\{\sigma_A,\tau_A\}$ is
a generating set of $M$. 

The new generators satisfy many relations that follow from
 the relations for $c_A$
and $e_A$. We notice some of them for future use.
 We use symbol
$\rho_A$ meaning either $\sigma_A$ or $\tau_A$.

(3.1) All $\sigma_A$ have degree 2 and thus
belong to the center of $M$. The elements $\tau_A$ are
degree 1 generators of an exterior subalgebra of $M$. In particular
$\tau_A^2=0$ and $\tau_A\tau_B=-\tau_B\tau_A$.
Thus all monomials in $\tau_A$ and $\sigma_A$ can be 
parametrized by the triples $(S,T,m)$ where $S,T\subseteq X$, a linear ordering
on $T$ is fixed and $m:S\to{\mathbf
Z_+}$. We put
$$\mu(S,T,m)=\prod_{B\in T}\tau_B\prod_{A\in S}\sigma_A^{m(A)}$$
where the first product is taken according to the order on $T$.
Notice that $\deg\mu(S,T,m)=|T|+2\sum_{A\in S}m(A)$.

(3.2) For every $A,B\in X$ we have
$$\tau_A\tau_B=\tau_A\tau_{A+B}-\tau_B\tau_{A+B}$$

\proof Represent $\tau_A=\Delta_{A,B}+\tau_{A+B}$ where $$\Delta_{A,B}=
\sum_{C\supseteq A,C\not\supseteq B}e_C.$$
Similarly $\tau_B=\Delta_{B,A}+\tau_{A+B}$. Since for every $C_1\supseteq A,
C_1\not\supseteq B$ and $C_2\supseteq B,C_2\not\supseteq
 A$ the set $\{C_1,C_2\}$ is not a
flag we have $e_{C_1}e_{C_2}=0$ whence $\Delta_{A,B}\Delta_{B,A}=0$.
Using besides that $\tau_{A+B}^2=0$ we obtain the result.

(3.3) For every $A\in X$ we have
$$\sigma_A^{\dim A}=0$$
and 

$$\rho_A\sigma_B^k=\rho_B\sigma_B^k$$
for every $A\subseteq B$ and $k\geq \dim B-\dim A$.

\proof The first relation follows immediately from (2.1) if one takes
$Y=\emptyset$ and $B=A$. In order to prove the second one it
 suffices to take
$k=\dim B-\dim A$.
Fix $A,B\in X$ such that $A\subset B$ and put $k=\dim B-\dim A$. We have
$$\rho_A\sigma_B^k=\sum_{C\supseteq A}a_C\sigma_B^k=\sum_{C\supseteq A}a_C
(\sum_{D\supseteq B}c_D)^k\eqno(3.4)$$
where $a_A=c_A$ or $e_A$ respectively for $\rho_A=\sigma_A$ or $\tau_A$.
Consider the sum $b_C=a_C
(\sum_{D\supseteq B}c_D)^k$ for some $C$ such that $C\supseteq A$ but $C\not
\supseteq B$. If there is no inclusion between $C$ and $D$ then $a_Cc_D=0$.
So it suffices to sum over $D$ such that $D\supseteq C$, i.e. $D\supseteq B+C$. 
In other words $$b_C=a_C(\sum_{D\supseteq B+C}c_D)^k.$$
Now we have
$$\dim(B+C)-\dim C=\dim B-\dim B\cap C\leq\dim B-\dim A=k.$$
 Due to (2.1) $b_C=0$ and the summation in (3.4) can be taken over $C$
such that $C\supseteq B$.
                  \qed
                                   
\medskip
Now we are ready to exhibit a basis of $M$ consisting
of monomials in $\sigma_A$ and $\tau_A$.
\begin{proposition}
\label{newbasis}
Suppose that $S,T\in X$, $S\cup T$ is a flag (by inclusion), the ordering on
$T$ coincides with the inclusion order, $m(A)<\dim A-\dim 
{A'}$ for every $A\in S$
and finally $|T|+2\sum_{A\in S}m(A)=r$ for a fixed nonnegative $r$.
Then the set of monomials $\mu(S,T,m)$ 
is a linear basis of $M_r$ 
(we call the monomials in this set {\it basic} monomials).
\end{proposition}
\proof
Comparing the set of basic monomials of degree $r$ with the basis of $M_r$ from
Section 2 we see that they have the same number of elements. Thus it suffices
to proof that every monomial $\lambda=
\lambda(S,T,m)$ in $c_A$ and $e_A$ with $m(A)<\dim A-\dim {A'}$ and
 $2\sum_{A\in S}m(A)+|T|=r$
is a linear combination of basic monomials of degree $r$. 
Apply downward induction on weight of $\lambda$.
Clearly the weight of any monomial of degree $r$ is limited from above by
$nr$ where $n=\dim V$ and the monomial of weight $nr$ is either
$c_V^{r/2}=\sigma_V^{r/2}$ if $r$ is even or $e_Vc_V^{(r-1)/2}=
\tau_V\sigma_V^{(r-1)/2}$ otherwise.
Since this monomial is 0 if $\lfloor r/2 \rfloor\geq n$ and basic otherwise
we have the base of the induction. 

Let $\lambda=\lambda(S,T,m)$ 
be an arbitrary monomial in $c_A$ and $e_A$ of
degree $r$ with $m(A)<\dim A-\dim {A'}$ and ${\rm wt}(\lambda)<nr$.
Consider the polynomial
$p=\lambda-\mu(S,T,m)$ in $c_A$ and $e_A$. Clearly every monomial of $p$ has
a weight larger than the weight of $\lambda$ whence the same is true for each
monomial of the decomposition of $p$ into the basis of Section 2.
Thus the induction works.
               \qed 

\bigskip
\section{Critical monomials}
\bigskip
Whenever we consider a monomial $\mu(S,T,m)$
in the rest of the paper, if $S\subset T$ we extend $m$ to $T$ by 0.
\begin{definition}
A basic monomial $\mu(S,T,m)$ is critical if $T\supseteq S$ and
\break
$m(A)=\dim A-\dim {A'}-1$ for every $A\in
T$ (so $S=\{A\in T|\dim A-\dim {A'}>1\}$). 
\end{definition}

\medskip
Notice that the set of all
critical monomials is linearly independent in $M$. Notice also that a
critical monomial is defined by a flag $T$ provided with the dimensions of
its elements.
 Denote by $c\mu(T)$
the critical monomial
defined by $T$. We have
$$\deg c\mu(T)=2\dim {A(T)}-|T|$$
where $A(T)$ is the maximal element of $T$. Denote by $CM$ the linear
subspace of $M$ spanned by all the critical monomials.

\begin{proposition}
\label{diff}
Let $T=\{A_1,A_2,\ldots,A_k\}$ where $A_i\subset A_{i+1}$ and let
$T_i=\{A_1,\ldots,\hat A_i,\ldots, A_k\}$ for every $i=1,2,\ldots,k$. Then
$$d(c\mu(T))=\sum_{i=1}^{k-1}(-1)^{i}c\mu(T_i).$$
\end{proposition}
\proof
To simplify notation put $\rho_i=\rho_{A_i}$
 and $m(i)=\dim {A_i}-\dim {A_{i-1}}-1$. 
Using the Leibniz formula for
$d$ and the relations (3.3) repeatedly we have
$$d(c\mu(T))=d(\tau_1\cdots\tau_k\sigma_1^{m(1)}\cdots\sigma_k^{m(k)})$$
$$=\sum_{i=1}^k(-1)^{i-1}
\tau_1\cdots\hat\tau_i\cdots\tau_k\sigma_1^{m(1)}\cdots
\sigma_i^{m(i)+1}\cdots\sigma_k^{m(k)}$$
$$=\sum_{i=2}^k(-1)^{i-1}\tau_1\cdots\widehat{\tau_{i-1}}
\cdots\tau_k\sigma_1^{m(1)}
\cdots\widehat{\sigma_{i-1}^{m(i-1)}}\sigma_i^{m(i)+m(i-1)+1}\cdots\sigma_k
^{m(k)}$$
$$=\sum_{i=1}^{k-1}(-1)^ic\mu(T_i).$$
     \qed

\medskip
\begin{corollary}
\label{invariant}
The linear space $CM$ is a subcomplex of $M$.
\end{corollary}

\medskip
\begin{remark}
\label{end}
Proposition \ref{diff} implies that the complex $CM$ affords the grading 
$CM=\oplus_{A\in X}CM_A$ where $CM_A$ is generated by monomials $c\mu(T)$
with $A(T)=A$. The deletion of $A$ from all the flags gives an isomorphism
of degree -1
from $CM_A$ 
to the usual flag complex $F(A)$ of the poset $(0,A)$
(see Section 6).
\end{remark}
Our next goal is to prove that $CM$ is quasi-isomorphic to $M$.
Denote by $CM^{\perp}$ the linear space spanned by all the basic monomials
that are not critical. Clearly $M=CM\oplus CM^{\perp}$ as a linear space.
The next lemma proves that it is true in the category of complexes also.

\begin{lemma}
\label{complement}
The space $CM^{\perp}$ is a subcomplex of $M$.
\end{lemma}
\proof
Let $\mu=\mu(S,T,m)\in CM^{\perp}$. This means that $\mu$ satisfies one of the
following conditions.

1. $S\not\subseteq T$, i.e., there exists $A\in X$ such that $\sigma_A$
divides $\mu$ and $\tau_A$ does not. Then every monomial in $d(\mu)$
satisfies the same condition whence $d(\mu)\in CM^{\perp}$.

2. $S\subseteq T$ and there exists $A\in T$ such that $m(A)<\dim A-\dim{A'}-1$. 
Suppose $T=\{A_1,\ldots,A_k\}$ and $A=A_i$. As
in the proof of Proposition \ref{diff} we put $\rho_{A_j}=\rho_j$, $m(A_j)=
m(j)$ and have
$$d(\mu)=\sum_{j=2}^k(-1)^{j-1}\mu_j$$ where 
$$\mu_j=\tau_1\cdots\hat\tau_j\cdots\tau_k\sigma_1^{m(1)}\cdots
\sigma_j^{m(j)+1}\cdots\sigma_k^{m(k)}.$$
Clearly $\mu_j=\mu(T_j,S^j,m^{(j)})$ for some $S^j$ and a function $m^{(j)}$.

Consider several cases. If $m(j)<\dim A_j-\dim A_{j-1}-1$ then $d^{(j)}(\mu)$ is basic
and $S^j\not\subseteq T_j$ whence $d^{(j)}(\mu)\in CM^{\perp}$. Suppose
$m(j)=\dim A_j-\dim A_{j-1}-1$. Using the same computation as in the proof of
Proposition \ref{diff} we have
$$\mu_j=
\tau_1\cdots\hat\tau_{j-1}\cdots\tau_k\sigma_1^{m(1)}\cdots
\widehat{\sigma_{j-1}^{m(j-1)}}\sigma_j^{m(j)+m(j-1)+1}\cdots\sigma_k^{m(k)}
.$$
The monomial $\mu_j$ is again basic. If $j\not=i+1$ then
the exponent of $\sigma_i$ in $\mu_j$ is $m(i)$ whence $\mu_j\in
CM^{\perp}$. Finally if $j=i+1$ then the exponent of $\sigma_j$ 
in $\mu_j$ is 
$$m(j)+m(i)+1<\dim A_j-\dim A_i-1+\dim A_i-\dim A_{i-1}-1+1$$
$$=
\dim A_j-\dim A_{j-2}-1$$
whence again $\mu_j\in CM^{\perp}$. That completes the proof.
          \qed

\medskip
Now we are to construct a homotopy on $CM^{\perp}$ between the identity map
and 0.

For every basic monomial $\mu=\mu(S,T,m)$ and $A\in S\setminus T$
put $h_A\mu=
\mu(S,T\cup\{A\},m)/\sigma_A$.
Then put $h\mu=\sum_{A\in S\setminus T}
(-1)^{(A,T)}
h_A\mu$ 
 where $(A,T)=|\{B\in T|B\subset A\}|$.
Also for every $B\in T$ put 
$d^B(\mu)=\mu(S,T\setminus\{B\},m)\sigma_B$.
Notice that $d(\mu)=\sum_{B\in
T}(-1)^{(B,T)}d^B(\mu)$ where $(B,T)=|\{C\in T|C\subset B\}|$.
Call $B\in T$ {\it critical for $\mu$}
if $m(B)=\dim B-\dim{B'}-1$ (as usual put $m(B)=0$ if $B\not\in S$).
Denote by $CT$ the set of all critical elements of $T$.
Now put $|\mu(S,T,m)|=|S\cup T\setminus CT|$. Notice that $|\mu|\not=0$ if
and only if $\mu\in CM^{\perp}$.

\begin{proposition}
\label{homotopy}
The linear map $CM^{\perp}\to CM^{\perp}$ defined by $\mu\mapsto (1/{|\mu|})
h\mu$ is a homotopy between the identity map and 0.
\end{proposition}
\proof
It suffices to check that for every basic monomial $\mu$ we have
$$hd(\mu)+d(h\mu)=|\mu|\mu.\eqno(4.1)$$
To make the check easier we 
state several simple
observations. In them and in the rest of the proof we always
mean that $\mu=\mu(S,T,m)$. Also all the equalities mean in particular that
the expressions in them are defined.

(a) For every $A\in S\setminus T$
and $B\in T\setminus CT$ we have $h_Ad^B(\mu)=d^B(h_A\mu)$.

(b) For every $A$ and $B$ as in (a) we have $d^A(h_A\mu)=\mu=h_Bd^B(\mu)$.

(c) For every $A\in S\setminus T$
the set of all critical for $h_A\mu$ elements coincides with $CT$.

(d) For $B\in CT$, let
 $\mu(S',T',m')$ be the basic monomial equal to $d^B(\mu)$. Then $S'\setminus
T'=S\setminus T$.

(e) Let again $B\in CT$ and $A\in S\setminus T$. Then $h_Ad^B(\mu)=
d^B(h_A\mu)$.

The statements (a)-(c) are straightforward. To check (d) and (e) 
apply relations of type (3.3) several times and notice that $T\setminus
T'=S\setminus S'$.

Combining (a)-(e) and using the rule for signs in the Leibniz property 
of $d$ one gets (4.1).
           \qed

\medskip
\begin{corollary}
\label{qi}
The embedding $CM\subset M$ is a quasi-isomorphism of the complexes.
\end{corollary}
\bigskip
\section{Multiplicative structure of $CM$}
\bigskip
In this section, we prove that $CM$ is a subalgebra of $M$ and 
compute the product of two critical monomials explicitly in terms of the
poset $X$ labeled by dimension.

In fact, we do first a more general computation. Consider the ${\bf C}$-linear
space $W$ 
of all formal linear combinations of the pairs $(T,m)$ where $T\subseteq X$ is a
flag and $m:T \to {\mathbf N}$. We want to convert $W$ into an associative
algebra. For that we recall the shuffle product on the free Abelian group
generated by flags (cf. \cite{OT}, 3.4).
 Let $T_1=(A_1\subset
A_2\subset\cdots\subset A_p)$ and $T_2=(B_1\subset B_2\subset\cdots\subset
B_q)$ and $m_i:T_i\to{\mathbf N}$ ($i=1,2$). Denote by $\lambda$ 
the natural operator making a flag from an arbitrarily linearly
ordered subset of $X$.
More precisely
$$\lambda(C_1,C_2,\ldots,C_k)=(C_1,C_1+C_2,\ldots,C_1+C_2+\cdots+C_k)$$
if all the elements in the right-hand side are distinct and 0 (of the group)
otherwise.
Then define
$$T_1\circ T_2=\sum_{\pi}({\rm sign}\pi)\lambda\pi(T_1\cup T_2)$$
where $\pi$ runs through all $(p,q)$-shuffles of $T_1\cup T_2$ (these are the
permutations preserving the order among $A_i$'s and the order among $B_j$'s).
It is easy to check that the defined multiplication is associative (see
\cite{OT}, p.88).
For  the future use put for a $(p,q)$-shuffle $\pi$ of $T_1\cup T_2$ 
$$(T_1\cup T_2)^{\pi}=\lambda\pi(T_1\cup T_2).$$

We want to 
extend this multiplication to $W$. For that fix a $(p,q)$-shuffle $\pi$
and notice that the flag $(T_1\cup T_2)^{\pi}$ consists of some of
the spaces $A_i+B_j$
where $0\leq i\leq p$ and $0\leq j\leq q$ (we always put $A_0=B_0=0$). 
Notice also that the predecessor of $A_i+B_j$
is either $A_{i-1}+B_j$ or $A_i+B_{j-1}$ (where the predecessor of
the smallest element $A_1=A_1+B_0$ or $B_1=A_0+B_1$ is taken to be $A_0+B_0=0$).
Define $(m_1\circ m_2)^{\pi}:(T_1\cup T_2)^{\pi}\to{\bf N}$ via
$$(m_1\circ m_2)^{\pi}(A_i+B_j)=m_1(A_i)\ {\rm or}\ m_2(B_j)$$
if the predecessor of $A_i+B_j$ is $A_{i-1}+B_j$ or $A_i+B_{j-1}$ respectively.
Now put
$$(T_1,m_1)\circ (T_2,m_2)=\sum_{\pi}{\rm sign}(\pi)
((T_1\cup T_2)^{\pi},(m_1\circ m_2)^{\pi}
)$$
where $\pi$ again runs through the group of $(p,q)$-shuffles of $T_1\cup T_2$.
The associativity of this product follows easily from the associativity of
the shuffle product of flags. Hence this product defines a structure of an
associative algebra on $W$ (with the identity represented by the empty flag).

Define the linear map $f:W\to M$ via $f(\emptyset)=1$ and
$$f((T,m))=\tau_{A_1}\cdots\tau_{A_p}\prod_{i=1}^p\sigma_{A_i}^{m(A_i)}$$
for $T=(A_1\subset\cdots\subset A_p)$.

The main theorem of this section is as follows.

\begin{theorem}
\label{product}
The map $f$ is a homomorphism of the algebras.
\end{theorem}

\proof Fix $(T_1,m_1)$ and $(T_2,m_2)$ from $W$ where
$T_1=(A_1\subset\cdots\subset A_p)$ and $T_2=(B_1\subset\cdots\subset B_q)$.
We need to prove that
 $$f(T_1,m_1)f(T_2,m_2)=f((T_1,m_1)\circ (T_2,m_2)).\eqno(5.1)$$
First we consider particular cases.

(a) Suppose $m_i=0$ ($i=1,2$)
 and $p=q=1$. Then (5.1) is a relation of type (3.2).

(b) Suppose $m_i=0$ ($i=1,2$)
 and $p=1$. Putting $A=A_1$ we have (5.1) in the form
$$\tau_A\tau_{B_1}\cdots\tau_{B_q}=\sum_{r=0}^k(-1)^r
\tau_{B_1}\cdots\tau_{B_r}\tau_
{A+B_r}\cdots\tau_{A+B_q}.\eqno(5.2)$$
We prove (5.2) by induction on $q$ using (a) as the base for $q=1$. For $q>1$ 
by the inductive hypothesis we have
$$\tau_A\tau_{B_1}\cdots\tau_{B_q}=\sum_{r=0}^{q-1}(-1)^r
\tau_{B_1}\cdots\tau_{B_r}
\tau_{A+B_r}\cdots\tau_{A+B_{q-1}}\tau_{B_q}.\eqno(5.3)$$
Applying relation (3.2) to the last two factors of (5.3) we can
substitute instead of them 
$$\tau_{A+B_{q-1}}\tau_{A+B_q}-\tau_{B_q}\tau_{A+B_q}.\eqno(5.4)$$
While substituting the second summand of (5.4) we get $q-1$ monomials 
ending with
$$\tau_{A+B_{q-2}}\tau_{B_q}\tau_{A+B_q}\eqno(5.5)$$
and the monomial
$$\tau_{B_1}\cdots\tau_{B_{q-1}}\tau_{B_q}\tau_{A+B_q}\eqno(5.6)$$
with coefficient $(-1)^q$.
Applying again relation (3.2) to the first two factors from (5.5) and using
that $A+B_{q-2}+B_q=A+B_q$ and $\tau_{A+B_q}^2=0$ we obtain that all $q-1$
monomials of the first kind are 0. Knowing that and substituting (5.6) and
(5.4) into (5.3) we obtain (5.2).

(c) Suppose $m_i=0\ (i=1,2)$.
 The proof is straightforward by induction on $p$ using (b)
as the base for $p=1$ and the associativity of the shuffle product.

(d) To finish the proof we need the following claim:
$$\tau_A\tau_{A+B}\sigma_{B}=\tau_A\tau_{A+B}\sigma_{A+B}$$
for every $A,B\in X$. The proof is similar to the proof of relation
(3.2) using also that $\tau_{A+B}^2=0$. The details are left to the reader.

(e) Consider the general case. By (c) the left hand side of (5.1) can be
rewritten as 
$$\sum_{\pi}({\rm sign}\pi)
f((T_1\circ T_2)^{\pi},0)\prod_{i=1}^p\sigma_{A_i}^{m_1(A_i)}\prod_{j=1}^q
\sigma_{B_j}^{m_2(B_j)}.\eqno(5.7)$$
Fix a $\pi$ and denote the respective monomial of (5.7) by $s(\pi)$.
Suppose $A_i+B_j\in (T_1\circ T_2)^{\pi}$. If $A_{i-1}+B_j\in (T_1\circ T_2)
^{\pi}$ (and thus precedes $A_i+B_j$) then apply (d) to 
$$\tau_{A_{i-1}+B_j}\tau_{A_i+B_j}\sigma_{A_i}^{m_1(A_i)}\eqno(5.8)$$
$m_1(A_i)$ times.
Since $A_{i-1}+B_j+A_i=A_i+B_j$ we can substitute
$$\tau_{A_{i-1}+B_j}\tau_{A_i+B_j}\sigma_{A_i+B_j}^{m_1(A_i)}$$
instead of (5.8). The case where $A_i+B_{j-1}\in (T_1\circ T_2)
^{\pi}$ can be handled similarly. 
Since eventually every factor of the form $\sigma_C^k$ 
is used we obtain
$$s(\pi)=f((T_1\circ T_2)^{\pi},(m_1\circ m_2)^{\pi})$$
that completes the proof.
     \qed

\medskip
\begin{corollary}
\label{cproduct}
For every two flags $T_1$ and $T_2$ with the maximal elements $A$ and $B$
respectively we have
$$c\mu(T_1)c\mu(T_2)=\sum_{\pi}({\rm sign}\pi)c\mu(T_1\circ T_2)^{\pi}$$
if $A\cap B=0$ (equivalently $\dim (A+B)=\dim A+\dim B$) and
$$c\mu(T_1)c\mu(T_2)=0$$
otherwise. Here again $\pi$ runs through all the $(|T_1|,|T_2|)$-shuffles
of $T_1\cup T_2$.
\end{corollary}

\proof Fix $T_1$ and $T_2$ as in the proof of the theorem above and fix a
$(p,q)$-shuffle $\pi$. Define $m_k:T_k\to{\ints_+}$ ($k=1,2$) via
$m_1(A_i)=\dim A_i-\dim A_{i-1}-1$ and $m_2(B_j)=\dim B_j-\dim B_{j-1}-1$. 
By definition
 $c\mu(T_k)=f((T_k,m_k))$ ($k=1,2$).

Suppose first that $A_p\cap B_q=0$. Let $A_i+B_j\in
(T_1\circ T_2)^{\pi}$. If its predecessor is $A_{i-1}+B_j$ then we have
$$(m_1\circ m_2)^{\pi}(A_i+B_j)=m_1(A_i)=\dim A_i-\dim A_{i-1}-1$$
$$=
\dim(A_i+B_j)-\dim(A_{i-1}+B_j)-1\eqno(5.9)$$
since $A_i\cap B_j=A_{i-1}\cap B_j=0$. The case where the predecessor
is $A_i+B_{j-1}$ can be handled similarly. Thus $$f((T_1\circ
T_2)^{\pi},(m_1\circ m_2)^{\pi})=c\mu((T_1\circ T_2)^{\pi})$$ that completes
the proof of the first case.

Now suppose that $A_p\cap B_q\not=0$. Then there exists $A_i+B_j\in (T_1\circ
T_2)^{\pi}$ that is the smallest (by inclusion) element such that $\dim
(A_i+B_j)<\dim A_i+\dim B_j$. To simplify notation put $m=(m_1\circ
m_2)^{\pi}$, $T=(T_1\circ T_2)^{\pi}$ and denote all the elements of $T$
from the
smallest to $A_i+B_j$ by $C_1,C_2,\ldots, C_k$. Then similar calculations as
in (5.9) give $m(C_i)=\dim C_i-\dim C_{i-1}-1$ for $i<k$ and $m(C_k)\geq \dim
C_k-\dim C_{k-1}$. Applying several times relation (3.3) we see that
$f((T,m))$ is equal to a monomial having $\sigma_{C_k}^{\ell}$ in it with
$\ell\geq\dim C_k$. Thus $f((T,m))=0$, that completes the proof.
               \qed

\medskip
Now Corollaries \ref{cproduct} and \ref{qi} imply the following.

\begin{corollary}
\label{subalgebra}
$CM$ is a (differential graded) subalgebra of $M$ and
a rational model of $C(X)$.
\end{corollary}

\medskip
\begin{remark}
\label{ints}
DGA $CM$ has a natural integer version (and a differential graded subring)
$CM(\ints)$ generated as a free Abelian group by the monomials
$c\mu(T)$.
\end{remark}

\bigskip
\section{Algebra $H^*(CM)$}
\bigskip

In this section, we use 
the combinatorial description of algebra $CM$ from the previous section to
give a
combinatorial description of algebra $H^*(CM)$.

In the rest of the paper, we will work much with homology of lattices. Let us
introduce the notation that we will use. A lattice $L$ has two
operations: the least upper bound $\vee$ (join)
 and the greatest lower bound $\we$ (meet).
For any subset $\sigma\subseteq L$ put
$$\bigvee(\sigma)=\bigvee_{A\in\sigma}A.$$
The smallest and largest elements of a lattice will be denoted by 0 and 1
respectively unless more specific symbols are appropriate. The set of all the
atoms of $L$ will be denoted by $\A=\A(L)$. An arbitrary linear ordering is
fixed on this set and all subsets of atoms will be ordered by the induced
order.  For every $T\in
L$ the set of complements of $T$ is $\C(T)=\{S\in L|S\vee T=1,S\wedge T=0\}$.
The homology of an
arbitrary poset $P$ is the homology of its complex $F(P)$
of flags. For a lattice $L$,
it is customary to write $H_*(L)$ for the homology of the poset
$L_0=
L\setminus\{0,1\}$. This homology can be computed also as the homology of the
atomic complex $\Delta(L)$ 
on $\A$ whose simplexes are the subsets $\sigma\subseteq \A$ with
 $\bigvee(\sigma)<1$.
In fact $\Delta(L)$ is homotopy equivalent to $F(L_0)$ (see \cite{BWa} and
Lemma \ref{fromatoms} below).
For any simplicial complex $\Delta$ we denote by $C(\Delta)$ its
chain complex over $\ints$.

We need to fix a particular homotopy equivalence of the chain
complexes $C(\Delta(L))$ and $C(F(L_0))$.
For that denote by $\Delta'(L)$ the barycentric subdivision of $\Delta(L)$
and identify $\Delta'(L)$ with $F(Q)$ where $Q$
is the poset of all nonempty simplexes of
$\Delta(L)$ ordered by inclusion. 
Denote by $\beta$ the standard homotopy equivalence $\beta:C(\Delta(L))\to
C(\Delta'(L))$ defined by 
$$\beta(\{A_1,\ldots,A_p\})=
\sum_{\delta}(-1)^{{\rm sign}
\delta}(\{A_{\delta(1)}\},
\{A_{\delta(1)},A_{\delta(2)}\},\cdots,
\{A_{\delta(1)},\ldots,A_{\delta(p)}\})$$
where $\{A_1,\ldots,A_p\}$ is an arbitrary simplex of $\Delta(L)$ (ordered
according to the fixed order on $\A$)
and the summation is taken over all the permutations $\delta$ of rank $p$ .
Also define the order preserving map $\gamma:Q\to L_0$ via $\gamma(\sigma)=
\bigvee(\sigma)$ for $\sigma\in\Delta(L)$ and keep the same symbol for the
respective chain map $C(\Delta'(L))=C(F(Q))\to C(F(L_0))$.
The map $\gamma$ is a homotopy equivalence since for every $C\in L_0$ the
poset $\gamma^{-1}(C)$ has the unique maximal element $\{A\in\A(L)|A\leq C\}$
(e.g., see \cite{Qu}). Thus we can register the following lemma.

\begin{lemma}
\label{fromatoms}
For a lattice $L$ 
the group homomorphism 
$f_L=\gamma\beta:C(\Delta(L))\to
C(F(L_0))$ given on atomic simplexes by
$$f(A_1,\ldots,A_p)=\sum_{\delta}(-1)^{{\rm sign}
\delta}(A_{\delta(1)}\leq
A_{\delta(1)}\vee A_{\delta(2)}\leq\cdots\leq
\bigvee_{i=1}^pA_{\delta(i)}),$$
where the summation is taken over all the permutations $\delta$ of rank $p$ and
every flag with repetitions considered to be
0, is a chain homotopy equivalence.
\end{lemma}

\medskip
We will also write $f_{L_0}$ for $f_L$.

Let us introduce several more pieces of notation.
For every $A,B\in L$ such that $B\leq A$ there is the lattice (closed
interval) $[B,A]=\{C\in
L|B\leq C\leq A\}$ whose homology is denoted simply by $H_*(B,A)$. 
The set of atoms of $[0,A]$ we denote by $\A(A)$, the flag complex of
$[0,A]_0=(0,A)$ by $F(A)$
and the atomic complex of $[0,A]$ by
$\Delta(A)$. 
 For any simplicial complex $\Delta$ we denote 
by $\SS(\Delta)$ its suspension. For any simplicial complexes $\Delta_1$ and
$\Delta_2$ we denote by $\Delta_1*\Delta_2$ their join. Also we use 
$\overline{\Delta}$ for the simplex on all the vertices of $\Delta$.
 For every $A,B\in L$ we 
denote by $\Delta(A,B)$
the subcomplex of $\Delta(A\vee
B)$ on the set $\A(A,B)=
\A(A)\cup \A(B)$. The joins of all subsets of $\A(A,B)$
form a sublattice  of $[0,A\vee B]$ that we denote by $L_{A,B}$
and we write $i(A,B)$ for the inclusion $L_{A,B}\subseteq [0,A\vee B]$.
We denote by $F_{A,B}$
the flag complex of $L_{A,B}$. 

Now we return to the poset $X$ of subspaces. When we apply to $X$
 notation and results
from lattice theory we always mean that 0 is adjoined to it.

\begin{lemma}
\label{join}
If $A,B\in X$ are such that $A\cap B=0$ then
$$\Delta(A,B)=(\overline{\Delta(A)}*\Delta(B))\cup (\Delta(A)*\overline{
\Delta(B)}).$$
\end{lemma}
\proof
Condition $A\cap B=0$ implies the similar condition on any
subspaces of $A$ and $B$. In particular, for $A_1\leq A$ and $B_1\leq B$ we
have $A_1\vee B_1=A\vee B$ only if $A_1=A$ and $B_1=B$. The result follows
immediately.                     \qed

\medskip
Lemma \ref{join} implies in particular that under the condition of the lemma,
$\Delta(A,B)$ is homotopy equivalent to $\SS(\Delta(A)*\Delta(B))$. To fix a
chain homotopy equivalence it is convenient to pass to $\Delta(A,B)'$. Denote
the vertices of the suspension by $a$ and $b$ and consider the chain map
$$\phi=\phi(A,B):
C(\SS(\Delta(A)*\Delta(B)))\to C(\Delta(A,B)')$$
defined by applying $\beta$ to simplexes from $\Delta(A)*\Delta(B)$
and mapping $a$ and $b$ to the barycenters of $\overline{\Delta(A)}$
and $\overline{\Delta(B)}$ respectively. For this map to be well-defined we
need to specify that the vertices of the suspension are taken the last among
vertices of simplexes containing them.
For instance, if $\{A_1\}\in\Delta(A)$, $\{B_1\}\in\Delta(B)$ and $A_1<B_1$ in
the ordering on atoms then
$$\phi(\{A_1,B_1,a\})$$
$$=(A_1,\{A_1,B_1\},D)-
(A_1,\overline{\Delta(A)},D)+(B_1,\{A_1,B_1\},D)$$
where 
$D=\overline{\Delta(A)}\cup\{B_1\}$.
Clearly $\phi$ is a chain homotopy equivalence.

Finally we need the following lemma well-known for arbitrary simplicial
complexes (e.g., see \cite{Mu}).

\begin{lemma}
\label{tensor}
For every $A,B\in X$ the group homomorphism 
$$\psi=\psi(A,B):
C(\Delta(A))\otimes C(\Delta(B))\to C(\SS(\Delta(A)*\Delta(B))),$$
given by 
$$\sigma\otimes\tau\mapsto (-1)^{\dim\tau+{\rm sign}
\epsilon(\sigma,\tau)} (\sigma\cup\tau\cup\{b\}-
\sigma\cup\tau\cup\{a\})$$
for $\sigma\in\Delta(A)$ and $\tau\in\Delta(B)$,
is a chain homotopy equivalence (of degree +2).
Here $\epsilon(\sigma,\tau)$ is the shuffle of $\sigma\cup\tau$ putting
all the elements of $\tau$ after elements of $\sigma$. 
\end{lemma}

\medskip
The sign for
$\psi(\sigma\otimes\tau)$ in Lemma \ref{tensor}
 is chosen so that the diagram in the proof of
Theorem \ref{cohomology} will be commutative.

If we put $\alpha(A,B)=\alpha=\phi(A,B)\psi(A,B)$ then
Lemmas \ref{join} and \ref{tensor} and the K\"unneth formula 
imply the following.

\begin{corollary}
\label{iso}
The map
$$\alpha_*: \tilde H_r((0,A),{\rats})\otimes \tilde H_s((0,B),{\rats})\to 
\tilde H_{r+s+2}(X_{A,B},{\rats})$$
is an isomorphism for every $r$ and $s$.
If the integer homology groups are free then the same is true over $\ints$.
\end{corollary}
 
\medskip

Now we are ready to describe  $H^*(C(X),\rats)=H^*(CM)$.
\begin{theorem}
\label{cohomology}
(i) For every $p$ we have isomorphisms
$$\tilde H^p(C(X),{\rats})\approx 
\tilde H^p(CM)\approx \oplus_{A\in X}\tilde H_{2\dim A-p-2}((0,A),
{\rats}) .$$
(ii) Under the isomorphism from (i) 
the multiplication on $H^*(C(X),{\rats})$ is given by 
the compositions

$$\tilde H_r((0,A),{\rats})\otimes \tilde H_s((0,B),{\rats})
\def\mapright#1{\smash{
   \mathop{\longrightarrow}\limits^{#1}}}
\mapright{\alpha_*} \tilde H_{r+s+2}(X_{A,B},\rats)\mapright{i_*}
\tilde H_{r+s+2}((0,A\vee
B),{\rats})$$
 if $A\cap B=0$ (equivalently $\dim(A\vee B)=\dim A+\dim B$)
 and 0 otherwise.
\end{theorem}

\proof
(i) follows immediately from the results of Chapter 4, in particular Remark
\ref{end} (and
recovers some results of \cite{GMP}).
Thus we
need to prove (ii) only.

In the rest of the proof all homology groups have the rational coefficients
and we suppress ${\rats}$. We know that the multiplication on $H_*(0,A)
\otimes H_*(0,B)$ is 0 unless $A\cap B=0$. Hence let us fix $A,B\in X$ such
that $A\cap B=0$.
 Then consider the following diagram

\smallskip
$$\def\mapright#1{\smash{
   \mathop{\longrightarrow}\limits^{#1}}}
  \def\mapleft#1{\smash{
   \mathop{\longleftarrow}\limits^{#1}}}
  \def\mapdown#1{\Big\downarrow
    \rlap{$\vcenter{\hbox{$#1$}}$}}
\matrix{C(\Delta(A))\otimes C(\Delta(B))&\mapright
{\psi(A,B)}&C(\SS(\Delta(A)*
\Delta(B)))&\mapright{\phi(A,B)}&C(\Delta(A,B)')\cr
\mapdown{f_A\otimes f_B}&&&&\mapdown{\gamma(A,B)}\cr
C(F_A)\otimes C(F_B)&&\mapright{\nu(A,B)}&&C(F_{A,B})\cr}$$

\smallskip
where $\nu(A,B)$ is induced by the shuffle
multiplication of flags from the previous section and $\gamma(A,B)$
is the natural map (introduced before Lemma \ref{fromatoms}) for the lattice
$X_{A,B}$. Notice that before one applies the shuffle multiplication to two
flags from $(0,A)$ and $(0,B)$ one should augment them by $A$ and $B$
respectively and then delete $A\vee B$ from all flags in the support of the
product.
A straightforward check shows that the diagram commutes.
Since the vertical maps are homotopy
equivalences,
passing to the homology completes the proof.
             \qed

\medskip
The statement (ii) of Theorem \ref{cohomology} makes sense over $\ints$.
I make the following conjecture.

\begin{conjecture}
\label{int}
Theorem \ref{cohomology} (ii) holds for the ring $H^*(C(X),\ints)$ (using
integer homology of the posets).
In particular this ring
is defined by the poset structure on $X$ 
(the ``combinatorics'' of the arrangement) 
and the
dimensions of its elements.
\end{conjecture}
\bigskip
\section{Geometric lattices}
\bigskip
In this section, 
we suppose that  $X$ (augmented by 0) is a geometric lattice,
i.e., its every element is the join of a set of atoms, for every $A\in
X$ all maximal flags ending on $A$ have the same length denoted by
${\rm rk}A+1$ and
$\rk(A\vee B)+\rk (A\wedge B)\leq \rk A+\rk B$ for every $A,B\in X$.
In this case we are able to give a more explicit description of $H^*(C(X))$
(all cohomology and homology in this and the following section are over
$\rats$).

By Folkman's theorem \cite{Fo} 
 $\tilde H_p(0,A)\not=0$ only for
$p=\rk A-2$.

The following observation will be helpful. 

\begin{lemma}
\label{rank}
For every $A,B\in X$ such that $A\cap B=0$ we have $\rk (A\vee B)=\rk A+\rk
B$.
\end{lemma}

\proof
Put $k=\rk A$, $\ell=\rk B$ and fix two flags in $X$
$$0=A_0<A_1<A_2<\cdots<A_k=A, \ 0=B_0<B_1<B_2<\cdots<B_{\ell}=B.$$
For any $i$, $0\leq i\leq \ell$, we have $A\cap B_i\subseteq A\cap B=0$
whence $\dim (A+ B_i)=\dim A+\dim B_i$. Thus we have the flag in $X$
$$0<A_1<A_2<\cdots <A_{k}=A<A+B_1<A+B_2<\cdots<A+B_{\ell}=A+B$$
that implies $\rk (A\vee B)\geq \rk A+\rk B$. The opposite inequality follows
immediately from the inequality for rank in geometric lattices.        \qed

\medskip
Thus the multiplication $\nu(A,B)$, which is 0 if $A\cap B\not=0$, is a map
$$\nu(A,B):\tilde H_{{\rm rk} A-2}(0,A)\otimes \tilde H_{{\rm rk} B-2}(0,B)\to
\tilde H_{{\rm rk}(A\vee B)-2}(0,A\vee B)\eqno(7.1)$$
otherwise.

To describe multiplication (7.1) more precisely and
 give a presentation of $H^*(C(X))$ we use the well known basis of
the
homology of geometric lattices (e.g., see \cite{OT}, p.144). 
 Recall
that a set $\sigma$ of atoms of a geometric lattice
 is called independent if $\rk \bigvee(\sigma)=|\sigma|$.
Clearly every subset of an independent set is again independent.
 Suppose $A\in X$ with $\rk A=p$. 
Fix an
arbitrary linear order on the set of atoms $\A(X)$.
Then for any independent 
subset
$\sigma=\{A_1,\ldots,A_p\}\subseteq\A(A)$ 
the cycle $\partial(\sigma)=\sum_{i=1}^p(-1)^i\sigma\setminus\{A_i\}$
 of $\Delta(A)$ defines a
non-zero element $\zeta_{\sigma}\in\tilde H_{p-2}(0,A)$
and all these elements generate the space. If a set $\sigma$ is dependent we will
use the symbol $\zeta_{\sigma}$ sometimes too but its value will be always 0.
 In general the
elements $\zeta_{\sigma}$ are not linearly
independent over $\ints$ even for independent sets $\sigma$. Linear relations
among them are generated by 
$$\sum_{j=1}^{p+1}(-1)^j\zeta_{\tau\setminus\{B_j\}}\eqno(7.2)$$
for every dependent set
$\tau=\{B_1,\ldots,B_{p+1}\}\subseteq\A(A)$.
Lemma \ref{rank} implies
 that if $A\cap B=0$ and $\sigma\subseteq \A(A)$ and $\tau
\subseteq\A(B)$ are
independent sets then the set $\sigma\cup \tau\subseteq \A(A\vee B)$ is again
independent.

Now Theorem \ref{cohomology} implies the following result.

\begin{proposition}
\label{geometric}
Suppose $X$ is a geometric lattice. Then the algebra \hfill\break
$H^*(C(X))$ is generated by
$\zeta_{\sigma}\in H^q(C(X))$ 
where $\sigma$ runs through all the independent sets of
atoms of $X$ and $q=2\dim\bigvee(\sigma)-|\sigma|$.
 Besides the linear relations (7.2), only other generating relations 
are
$$\zeta_{\sigma}\zeta_{\tau}=(-1)^{{\rm sign}\epsilon(\sigma,\tau)}
\zeta_{\sigma\cup\tau}$$
if $\bigvee(\sigma)\cap\bigvee(\tau)=0$ (where $\epsilon(\sigma,\tau)$ 
is the permutation of
$\sigma\cup\tau$ putting all the
elements of $\tau$ after elements of $\sigma$)
and
$$\zeta_{\sigma}\zeta_{\tau}=0$$
otherwise.
\end{proposition}

\medskip
\begin{remark}
\label{relations}
(i) In general, relations (7.2) are not linearly
independent.
One
can easily construct a linear resolution of the space of 
generators $\zeta_{\sigma}$ and
recover Folkman's formula for the dimension of $H(0,A)$ from \cite{Fo}.

(ii) Proposition \ref{geometric} shows in particular that 
there exists a quasi-iso-\hfill\break morphism
$H^*(C(X))\to M$ (of graded differential algebras with $d=0$ on 
\hfill\break
$H^*(C(X))$) whence
the space $C(X)$ is formal in this case.

(iii) In the case when $X$ is not only geometric but Boolean all the spaces
in (7.1) are 1-dimensional and the
multiplication is clearly an isomorphism
(cf. \cite{Fe}).
In general the multiplication is not an isomorphism.
\end{remark}

\bigskip
\section{Cohomology rings of ``$k$-equal'' manifolds} 
\bigskip
In this section, for each pair $(n,k)$ of integers
with $2\leq k\leq n$ and for each 
 $k$-subset $\omega\subseteq \{1,2,\ldots,n\}$
we consider the
subspace $V(\omega)\subset
{\bf C}^n=\{(x_1,\ldots,x_n)\}$ 
defined by the system of
equations $\sum_{i\in \omega}x_i=0$ and $x_j=0\ {\rm for}\ j\not\in\omega$.
The set of all $V(\omega)$ (for fixed $n$ and $k$) is denoted by $X_{n,k}$.
Notice that the annihilator $V(\omega)^0$ of $V(\omega)$ is
defined in $({\bf C}^n)^*=\{(z_1,\ldots,z_n)\}$
by the equations
$z_i=z_j$ for $i,j\in\omega$ that justifies the name $k$-equal manifold
introduced in \cite{BW} for
$M_{n,k}=C(X_{n,k})$. 

The linear structure of $H^*(M_{n,k})$ has been
 studied in detail in \cite{BW}.
Completing $X_{n,k}$ by the sums of its elements one obtains the lattice 
 $\Pi_{n,k}$ that is (naturally isomorphic to)
the lattice of all partitions of the set
$\{1,\ldots,n\}$
 whose blocks consist of either 1 element (trivial blocks)
or at least $k$ elements (non-trivial blocks).
In the rest of the paper, we will regard elements of $\Pi_{n,k}$ both as
subspaces and partitions.
Theorem 1.5 of \cite{BW} implies that $\tilde H_p(\Pi_{n,k})\not=0$
 if and only if
$p=n-3-t(k-2)$ for some integer $t$ such that $1\leq t\leq \lfloor
n/k\rfloor$. 

Although the lattice $\Pi_{n,k}$ is not geometric in general, it has a
similar property; there exist
generators of its homology groups represented each by a simplicial sphere
in the atomic complex. To prove that we need to use a recursive construction
 from \cite{BW}
hence to include into consideration more general lattices
introduced there.
For any non-negative integers $n,k$, and $\ell$ such that $n,k\geq 2$ and either
$k\leq n$ or $\ell>0$,
define $\Pi_{n,k}(\ell)$ as the family of all partitions of the set $\{1.\
\ldots,n\}$
whose every block $b$ satisfies at least one of the following requirements:

(i) $|b|=1$,

(ii) $|b|\geq k$,

(iii) $b\cap\{1,\ldots,\ell\}\not=\emptyset$.

\noindent
Ordered by refinement $\Pi_{n,k}(\ell)$ becomes a lattice and
$\Pi_{n,k}(0)=\Pi_{n,k}$.

First we need 
to recall the recursive construction
 used in \cite{BW} for the homology groups of
$\Pi_{n,k}(\ell)$. 
Suppose $n\geq 3$ and $k,\ell\leq n-1$.
 Let $T$ be a coatom of
$\Pi_{n,k}(\ell)$ having the non-trivial block $b=\{1,\ldots,n-1\}$.
Clearly there are two types of complements $S$ of $T$, both are atoms 
whose unique non-trivial blocks we denote by $b_1$.
Either $|b_1|=k$ and $b_1\cap\{1,\ldots,\ell\}=\emptyset$
 or $b_1=\{i,n\}$ where $i\leq\ell$.
  In the former case the poset
$[S,1]$ is isomorphic to $\Pi_{n-k+1,k}(\ell+1)$, in the latter case 
it is isomorphic to $\Pi_{n-1,k}(\ell)$. The subposet $(\Pi_{n,k}(\ell))_0=
(0,1)$ of $\Pi_{n,k}(\ell)$ is homotopy
equivalent to the wedge of suspensions of the posets $(S,1)$ where $S$ runs
over $\C(T)$ (see details and references in \cite{BW}).
In particular taking a basis of homology for each $[S,1]$ and applying
suspension one obtains a basis of homology of $\Pi_{n,k}(\ell)$. 
Notice that for each choice of $S$ the homotopy type of
the suspension of $(S,1)$ can be realized
as a subposet of $(0,1)$. For instance, one can take as such a
realization $(S,1)\cup\{S,T\}$ to which
at least one atom $A_0\leq T\wedge A$ for each $A\in (S,1)$ is added.
 We will call any of these
(homotopy equivalent to each other) posets the suspension of $(S,1)$ with
vertices $S$ and $T$. We will also use the same term for the suspension of
homology classes of $(S,1)$.

The result of the above 
recursive construction for $\Pi=\Pi_{n,k}$ can be described as follows. 
Fix 
a maximal flag $0=S_0<S_1<S_2<\cdots<S_{p+1}<1$ in $\Pi$ and a sequence of
elements $T_1,\ T_2,\ldots,T_{p}$ such that for every $i=1,2,\ldots,p$ 
the element  
$T_i$ has a non-trivial block of size $n-1$, 
$T_i\vee S_i=1$ and
$T_i\wedge S_i=S_{i-1}$.
A pair $(S',S_{p+1})$, where $S'$ is any atom of $[S_p,1]$ different from
$S_{p+1}$,
 defines an element $\alpha_p\in\tilde
H_0(S_p,1)$. Applying recursively suspensions with vertices $S_i$ and $T_i$
one obtains a sequence of $\alpha_i\in\tilde H_{p-i}(S_i,1)$. Then
elements $\alpha_0$ taken for all the above data generate $\tilde H_p(\Pi)$. 
The condition $\alpha_0\not=0$ puts restrictions on $S_p$ and $S_{p+1}$. In
fact, the variety of possibilities reduces to the following four cases (cf.
Lemma 4.4 from \cite{BW}).

(a) $n=k+1$ and $p=0$;

(b)  $S_p$ has only one non-trivial block of $n-k$ elements and $S_{p+1}$ has
besides the complementary block (of $k$ elements);

(c) $S_p$ has three blocks all together exactly one of which is non-trivial
 ($k\geq 3$);

(d) $S_p$ has three blocks all together at least two of which are
non-trivial.

Now we are ready to study the multiplication on local homology of
$\Pi$.
For every $p$ consider the homomorphism
$$\xi_p:\oplus_{\{A,B\}}\tilde H_p(\Pi_{A,B})\to\tilde H_p(\Pi)$$ 
induced by embeddings where $\{A,B\}$ runs through all the subsets of
$\Pi$ such that $A$ and $B$ are one-block partitions,
$A\vee B=1$ and $A\cap B=0$. The last two conditions can be
substituted by $b_A\cup b_B=\{1,\ldots,n\}$ and $|b_A\cap b_B|=1$ where $b_A$
 and $b_B$ are the
non-trivial blocks of $A$ and $B$ respectively.
Recall that the maximal $p$ with $\tilde H_p(\Pi)\not=0$ is $p=n-k-1$.
It is easy to see that $\xi_{n-k-1}=0$ since its domain is 0. 
On the other hand we have the
following.

\medskip
\begin{lemma}
\label{local}
For $p<n-k-1$ the homomorphism $\xi$ is surjective.
\end{lemma}

\proof
It suffices to fix generators of $\tilde H_p(\Pi)$ ($p<n-k-1$) and for each
$\alpha$ of them find $A$ and $B$ with required properties and such that
$\alpha\in{\rm Im}\xi$.

To fix an element $\alpha$
 of a generating set we fix sequences $(S_i)$ and $(T_i)$ as
above and, using suspension, construct recursively $\alpha_i\in\tilde H_{p-i}
(S_i,1)$ starting from $\alpha_p$ and finishing at $\alpha_0=\alpha$. 
The condition $p<n-k-1$ restricts possibilities
for $S_p$ to (b) and (d) only.
In case (b), we take as $S'$ any atom of $[S_p,1]$ different from $S_{p+1}$
($\alpha_p$ is essentially
unique in this case). In case (d), we are a little more
careful. Each $S_i$ has a block that contains the non-trivial block of
$S_1$. We call this block of $S_i$ distinguished. While choosing an atom $S'$
we avoid repetition of the distinguished block of $S_p$, i.e if this block
is a proper subset of another block in $S_{p+1}$ then in $S'$ 
this block is not a proper subset of any block.
Clearly such elements $\alpha_p$ generate $\tilde H_0(S_p,1)$.

Now we denote by $A_p$ the element of the pair $\{S_{p+1},S'\}$ that contains
the distinguished block of $S_p$ as a proper subset of 
another block and by $B_p$ the
other element of the pair. Our goal is to construct recursively
a sequence $(A_i,B_i)$
($i=p,p-1,\ldots,0$) with the following properties.

(i) $A_i,B_i>S_i$ and each has only one
non-trivial block ($b'_i$ and $b''_i$ respectively) with respect to $S_i$; 

(ii) $b_i'\cap b_i''$ is
a block of $S_i$; 

(iii) $A_i\vee B_i=1$; 

(iv) For $i>0$ only $A_i$ (and not $B_i$)
 has a block that contains the distinguished block
of $S_i$ as its proper subset;

(v) $\alpha_i$ is in the image of $\tilde H_i([S_i,1]_{A_i,B_i})$ in $\tilde
H_i(S_i,1)$.

Notice that $A_p$ and $B_p$ satisfy the conditions (i)-(v).
Now suppose that for some $i$, $0<i\leq p$,
$A_i$ and $B_i$ satisfying (i)-(v) are constructed.
The construction of $A_{i-1}$ and $B_{i-1}$ depends
 on $b_i'\cap b_i''$.
 Denote by $b_i$ the unique block
of $S_i$ that is non-trivial with respect to $S_{i-1}$ and let 
$b_i=\cup_jb_{ij}$
be its partition into blocks of $S_{i-1}$. 
If $b'_i\cap b''_i\not=b_i$ then define $A_{i-1}$ and $B_{i-1}$ 
as having the
non-trivial blocks $b'_{i-1}=b'_i$ and $b''_{i-1}=b''_i$ 
respectively, relative to $S_{i-1}$. If
$b'_i\cap b''_i=b_i$ then define $b''_{i-1}$ as in the previous case and put
 $b'_{i-1}=(b'_i\setminus
b_i)\cup b_{ij}$ for some $j$ such that $b_{ij}$ is a subset of the
non-trivial block of $T_i$.
The property (iv) of $A_i$ and $B_i$ guaranties that $A_{i-1}\in \Pi$.
In any case the conditions
(i)-(iv) for $i-1$ are satisfied by construction and we only need 
to check (v). 

To simplify notation put $P_i=[S_i,1]_{A_i,B_i}$ for every $i$
and denote the element not contained in the non-trivial block of $T_i$ by
$a_i$.
First notice that by construction either $A_{i-1}$ or $B_{i-1}\leq T_i$.
Suppose that $A_{i-1}\leq T_i$ since the other case can be handled
similarly (even more easily). This implies that $b_i\subset b''_{i-1}$.
Then let us prove that $T_i\in P_{i-1}$. 
Since $\Pi$ is atomic we have
$A_{i-1},B_{i-1}\in P_{i-1}$. Define $B'$ as the one-block
partition with the nontrivial block $b''_{i-1}\setminus\{a_i\}$.
Since $b_i\subset b''_{i-1}$ we have $B'\in \Pi$ whence $B'\in P_{i-1}$.
Since $A_{i-1}\vee B'=T_i$ we have $T_i\in P_{i-1}$. 

Our second observation is that $P_i\subset P_{i-1}$. To see that
it suffices to prove that any atom $C$ of $[S_i,1]$ such that $C\leq A_i$ or
$C\leq B_i$ is in $P_{i-1}$. Denote the only non-trivial with respect to $S_i$
block of $C$ by $b_C$. If $b_C\not\supset b_i$ then the atom $C_0$ of
$[S_{i-1},1]$ having only non-trivial block $b_C$ with respect to $S_{i-1}$
is in $P_{i-1}$ and
$C_0\vee S_i=C$. Suppose $b_C\supset b_i$ and $C\leq A_i$. Then $C_0=C\wedge
A_{i-1}$ is an atom of $P_{i-1}$ and again $C_0\vee S_i=C$. The case $C\leq
B_i$ is similar. In any case we have $C\in P_{i-1}$.

Summing up our above observations we see
that $P_{i-1}$ contains the suspension of $P_i$
with vertices $S_{i}$ and $T_i$ (cf.
 \cite{BWa} or Proposition 2.3  and p.290
of \cite{BW}).
Since by the inductive hypothesis (v) holds for $i$ we obtain it for $i-1$.

Putting $A=A_0$ and $B=B_0$ we complete the proof.
          \qed

\medskip
In order to exhibit more explicit generators of $H^*(M_{n,k})$ we use the
atomic complexes as in the previous section. We will always assume that the
atoms are linearly ordered and all sets of atoms are provided with the induced
ordering. 
By $\partial$ we denote the differential of the chain complex $C(L)$
of the simplicial
complex whose simplexes are all the sets of atoms.

\begin{definition}
A set $\sigma$ of atoms of a lattice $L$ is called independent 
if $\bigvee(\tau)\not=\bigvee(\sigma)$ for every 
proper subset $\tau$ of $\sigma$. It is called essential
if it is not a simplex in the atomic complex $\Delta(L)$, 
i.e. if $\bigvee(\sigma)=1$.
An element $\omega$ of $C(L)$ is called a relater in $L$
 if $\partial\omega$ can be represented as a linear combination of
independent essential sets and simplexes of the atomic complex, i.e., sets bounded
by elements less than 1.
\end{definition}

Any set $\sigma=\{A_1,A_2,\ldots,A_p\}$ of atoms of $L$ 
 defines the
 cycle $$\partial\sigma=
\sum_{i=1}^p(-1)^i\sigma\setminus\{A_i\}\in 
C_{p-2}(\Delta(L)).$$ Denote its homology class by
$\zeta_{\sigma}$.
If $\sigma$ is a simplex of $\Delta(L)$ then $\zeta_{\sigma}=0$.
Any relater
$\omega$ with $\partial\omega=\sum_ia_i\omega_i+\sum_jb_j\omega'_j$
 where $a_i,b_j\in\ints$, 
$\omega_i$ are independent essential sets and $\omega'_j$ are simplexes
defines the
linear relation $\sum_{i}a_i\zeta_{\omega_i}=0$.
This motivates the following. Denote by $G_p$ the linear space
generated by symbols $\eta_{\sigma}$ where $\sigma$ is either a simplex of
$\Delta(L)$ or
 an independent essential set of
atoms of $L$ and $|\sigma|=p+2$. Denote by $R_p$ the subspace of $G_p$
generated by the elements $\eta_{\sigma}$ where $\sigma$ is a simplex and by
$r_{\omega}=\sum_ia_i\eta_{\omega_i}$ where $\omega$ is a relater as above
 and $|\omega|=p+3$.

Fix $n,k,\ell$ such that $2\leq k\leq n$ and $\ell\geq 0$ and put
$\Pi=\Pi_{n,k}(\ell)$.

\begin{lemma}
\label{Claim}

(i) For every atom $S$ of $\Pi$ and a coatom $T$ with only one non-trivial
block such that $T\in\C(S)$
there exists an injective map $\theta:\A([S,1])\to \A(\Pi)\setminus\{S\}$ 
such that $\bigvee({\rm Im}\theta)\leq T$ and
$\theta(\tau)\cup\{S\}$ is independent
 (essential) in $\Pi$ if and only if
$\tau$ is independent (resp. essential) in $[S,1]$.
\end{lemma}

\proof

Let $S$ and $T$ be as in the statement of the lemma.
Consider the sets $W_A=\A(T\wedge A)$ for $A\in\A([S,1])$.
It is easy to see that for every $B\leq T\wedge A$ we have
$$B\vee S =A.$$
This implies that the sets $W_A$ are pairwise disjoint. Choosing arbitrarily
$\theta(A)\in W_A$ we obtain an injective map $\theta:\A([S,1])\to\A(\Pi)
\setminus\{S\}$. 

First let us prove the following property of $\theta$.
For every $\tau\subseteq\A([S,1])$
$$\bigvee(\theta(\tau)\cup\{S\})=\bigvee(\tau).\eqno(8.1)$$

Indeed put $\tau'=\theta(\tau\cup\{S\})$.
If $D\geq B$ for every $B\in\tau$ then $D>S$ and $D>
(T\wedge B)$ whence $D>A$
for every $A\in\tau'$. Conversely if $D\geq A$ for every $A\in\tau'$
then $D\geq(\theta(B)\vee S)=B$ for every $B\in\tau$.

Now it is easy to see that all the required properties of $\theta$ follow
from construction and (8.1).
           \qed

\medskip
\begin{remark}
\label{susp}
The construction of a map $\theta$ in Lemma \ref{Claim}
uses an atom $S$ and a coatom $T$. We will say that $\theta$ is constructed
with $S$ and $T$. Also it is clear from the construction
 that for every $\sigma\subseteq\A([S,1])$
the class $\zeta_{\theta(\sigma)}$ is the suspension of $\zeta_{\sigma}$ with
vertices $S$ and $T$.
\end{remark}

\medskip
Denote 
 the set of all essential independent subsets of atoms of a lattice 
$\Pi=\Pi_{n,k}(\ell)$ by
$\Io(\Pi)$.

\begin{proposition}
\label{property}
For every $p$, the linear map $\iota: G_p/R_p\to\tilde H_p(\Pi)$
 defined by $\iota(\eta_{\sigma})=\zeta_{\sigma}$
for all 
 $\sigma\in\Io(\Pi),\ |\sigma|=p+2$, 
is an isomorphism.
\end{proposition}

\proof
First we prove that $\iota$ is surjective.
We apply induction on $p$ (for all $n,k$ and
$\ell$ at the same time). As the base consider the case of $p=-1$. An essential
1-set $\sigma$ 
exists if and only if $(0,1)=\emptyset$ (i.e., $n=k,\ell=0$ or $n=2,\ell\geq
0$) and then it is unique and
independent. This is precisely the case when $\tilde H_{-1}(\Pi)\not=0$
 and $\zeta_{
\sigma}$
is a generator of this space.

Now suppose $p\geq 0$. If $k=n$ then the result is trivial. If $\ell\geq n$
or $n=2$
then $\Pi$ is a geometric lattice and the result follows from 
Proposition \ref{geometric}. 
Thus we can assume that $n\geq 3,\ \  k,\ell\leq n-1$ and apply
the above
recursive construction for 
$H_p(\Pi)$. 
Fix a coatom $T$ with only non-trivial block of $n-1$ elements
and atom $S$ such that
$S\in\C(T)$. Recall that $(S,1)$ is isomorphic to a poset
$\Pi_{n',k'}(\ell')$ for some value of the parameters $n',k'$ and $\ell'$. Thus
by the inductive hypothesis, one can find a basis of $\tilde
H_{p-1}(S,1)$ consisting of classes $\zeta_{\sigma}$ where $\sigma$ runs
through independent essential $p+1$-sets of atoms of $(S,1)$. 
Using Lemma \ref{Claim} (i) fix a map $\theta$ constructed with $S$ and $T$
and having property from that lemma. In particular if
$\sigma'=\theta(\sigma)$ then $\sigma'$ is independent and essential.
By Remark \ref{susp}
 $\zeta_{\sigma'}$ is the suspension of $\zeta_{\sigma}$ with the vertices
$S$ and $T$. Since
the suspensions of generators of all spaces $H_{p-1}(S,1)$ when $S$ runs
through
 $\C(T)$ form a generating set of $H_p(\Pi)$ the proof is complete.

Now we prove that $\iota$ is an isomorphism. For that consider the factor
complex $D=C(\Pi)/C(\Delta(\Pi))$. Notice that $D_p$ is generated by the
essential sets of atoms of cardinality $p+1$ and the differential on $D$ is
defined by
$$d(\sigma)=\sum(-1)^i\sigma\setminus\{A_i\}$$
where $\sigma=(A_1,\ldots,A_{p+1})$ and the summation is taken over all $i$
such that
$$\bigvee(\sigma\setminus\{A_i\})=1.$$
In particular an essential set is a cycle in $D$ if and only if it is
independent.
Applying the exact homological sequence of pair we see that the homomorphism
$\delta:\tilde H_{p+1}(D)\to\tilde H_p(C(\Delta))$ given by
$\delta(z)=\partial z$ for any cycle $z\in D_p$
 is an isomorphism. On the class $[\sigma])$ of an
independent essential set $\sigma$ we have
$\delta([\sigma])=\zeta_{\sigma}$. Thus the first part of the proof can be
restated as follows:
the classes of all the independent essential sets generate $H_*(D)$.
On the other hand, $R_p$ is naturally isomorphic to the intersection
$d(D_{p+2})\cap<\Io>$ in
$D_{p+1}$ and the result follows immediately.                  \qed

\medskip
The next lemma follows from Theorem \ref{cohomology} (ii) and \cite{BW},
p.296.

\begin{lemma}
\label{disconnected}
Suppose $U\in \Pi_{n,k}$ and has precisely $s$ non-trivial blocks\hfill\break
$b_1,\ldots,b_s$. Denote
by $U_i$ the element of $\Pi_{n,k}$ having only one non-trivial block
$b_i$. Then for every $p$ the repeated multiplication induces the 
isomorphism
$$\oplus_{p_1+\cdots
+p_s=p-2(s-1)}
\tilde H_{p_1}(0,U_1)\otimes\cdots\otimes\tilde 
H_{p_s}(0,U_s)\cong\tilde H_p(0,U).\eqno(8.2)$$
\end{lemma}

\medskip

\begin{definition}
\label{}
Let $\sigma$ be a set of one-block partitions (e.g. atoms) from $\Pi_{n,k}$.
Consider the collection $\D$ of all the subsets $\tau$ of $\sigma$ such that the
non-trivial blocks of elements of $\tau$ are pairwise disjoint.
Then the rank of $\sigma$ ($\rk(\sigma)$) is the maximal 
cardinality of elements of $\D$.
\end{definition}

\medskip
Here is the main result of the section.
\begin{theorem}
\label{k-equal}
(i) The space $H^*(M_{n,k})$ is generated by the classes \hfill\break
$\zeta_{\sigma}
\in H^{2\dim U-|\sigma|}(M_{n,k})$
where $\sigma$ runs through all the independent sets of atoms of 
$\Pi_{n,k}$ and
$U=\bigvee(\sigma)$.
The linear relations among $\zeta_{\sigma}$ 
 are generated by 
$r_{\omega}$ for the relaters $\omega$ of all $[0,A]$ $(A\in X)$.

(ii) The product $\zeta_{\sigma}\zeta_{\tau}=0$ if $\bigvee(\sigma)\cap
\bigvee(\tau)=0$ (equivalently
$\dim\bigvee(\sigma\cup\tau)\not=\dim\bigvee(\sigma)+\dim\bigvee(\tau)$) and
$\zeta_{\sigma}\zeta_{\tau}=(-1)^{{\rm sign}\epsilon(\sigma,\tau)}
\zeta_{\sigma\cup\tau}$ otherwise
where $\epsilon(\sigma,\tau)$ is as above the shuffle of
$\sigma\cup\tau$
 putting all the elements of $\tau$ after elements of $\sigma$.

(iii) The linear relations from (i) imply that $\zeta_{\sigma}=0$ unless
$$|\sigma|=n(\sigma)-\rk(\sigma)(k-2)-s,\eqno(8.3)$$
where $s$ is the number of non-trivial blocks in the partition
$\bigvee(\sigma)$ and $n(\sigma)$ is the number of points in all these blocks.
In particular $\zeta_{\sigma}\zeta_{\tau}=0$ unless
$\rk(\sigma\cup\tau)=\rk(\sigma)+\rk(\tau)$.

(iv) As an algebra, $H^*(M_{n,k})$ is generated by the classes $\zeta_{\sigma}$ for
independent sets of atoms of rank 1.
\end{theorem}

\proof
(i) and (ii) follow immediately from Proposition \ref{property} and Theorem
\ref{cohomology}.

(iii) Let $\sigma$ be a independent set of atoms, $r=\rk(\sigma)$,
$U=\bigvee(\sigma)$ and $U_1,\dots,U_s$ as in Lemma \ref{disconnected}.
Then $\sigma$ is a disjoint union of $\sigma_i$ ($i=1,\ldots,s$) where 
$\sigma_i$ is independent and $\bigvee(\sigma_i)=U_i$. Since
$r=\sum_{i=1}^s\rk(\sigma_i)$ and $\zeta_{\sigma}=\zeta_{\sigma_1}
\otimes\cdots\otimes\zeta_{\sigma_s}$ under identification (8.2) it suffices
to consider the case where $s=1$, i.e. $U$ has only one non-trivial block of
size $n(\sigma)$. It is easy to compute that
$$n(\sigma)\geq |\sigma|+r(k-2)+1. \eqno(8.4)$$

Now the recursive construction before Lemma \ref{fromatoms} shows that we can
assume without any loss of generality that
there exists a flag in the support of
$f_{[0,U]}(\sigma)$ that is maximal in $(0,U]$. 
If this maximal flag 
has an element with $t$ non-trivial blocks and no elements with $t+1$ ones
then its length is 
$$t+(t-1)+(n(\sigma)-tk)=n(\sigma)-t(k-2)-1$$
($t$ steps to
``create'' new blocks, $t-1$ steps to glue blocks to each other and
$n(\sigma)-tk$ steps to add a point to a block). Since for flags in the
support of $f_{[0,U]}(\sigma)$ the maximal possible $t$ is $r$ we have
$$|\sigma|\geq n(\sigma)-r(k-2)-1.\eqno(8.5)$$
Comparing (8.4) and (8.5) we obtain (8.3).

The second part of (iii) follows since the equality $\dim \bigvee(\sigma)+
\dim\bigvee(\tau)=\dim \bigvee(\sigma\cup\tau)$ implies $n(\sigma\cup\tau)=
n(\sigma)+n(\tau)-1$ for $\sigma$ and $\tau$ such that $\bigvee(\sigma)$ and
$\bigvee(\tau)$ have only one non-trivial block each.

(iv) Fix a independent set $\sigma$ of atoms of rank $r>1$ 
and let again $s$ be the number of non-trivial blocks
of $U=\bigvee(\sigma)$.
Using Lemma \ref{disconnected} again we can assume that $s=1$.
Suppose $\zeta_{\sigma}\not=0$.
Since $r>1$ we have from (iii) that
$\zeta_{\sigma}\in\tilde H_p(0,U)$ with
$$p=|\sigma|-2=n(\sigma)-r(k-2)-3<n(\sigma)-k-1.$$
Now due to Lemma \ref{local} and
Theorem \ref{cohomology} (ii), $\zeta_{\sigma}$ is a linear combination of
classes $\zeta_{\tau}\zeta_{\omega}=\pm\zeta_{\tau\cup\omega}$ 
where for each pair $(\tau,\omega)$
there exist $A$ and $B\in [0,U]$ such
that $A\vee B=U$, $A\cap B=0$,
$\zeta_{\tau}\in H_q(0,A)$ and $\zeta_{\omega}\in H_{p-q}(0,B)$ for some
$q$. Since $\bigvee(\tau\cup\omega)=U$ we obtain from (iii) that
$\rk (\tau\cup\omega)=r$.
Using (iii) again we have $\rk(\tau)<r$ and $\rk(\omega)<r$. Now downward
induction on $r$ 
completes the proof.                 \qed

\medskip
\begin{remark}
\label{rank1}
For any specified values of $n$ and $k$, the linear relations among
$\zeta_{\sigma}$ can be used to find a basis of $\tilde 
H_p(0,\bigvee(\sigma))$ ($p=|\sigma|-2$).
In the case $\rk(\sigma)=1$ a basis is easy to describe in general
that gives a generating set for the algebra $H^*(M_{n,k})$.

Indeed let the only non-trivial block of $U=\bigvee(\sigma)$ is $b=\{1,2,\ldots,
n(\sigma)\}$. Fix a subset $a$ of $b_0=b\setminus \{1\}$
 with $|a|=k-1$ and put $\sigma(a)=
\{A_i|i\in b\setminus a\}$ where $A_i$ is the atom whose only non-trivial
block is $a\cup\{i\}$. Notice that $|\sigma|=n(\sigma)-k+1$ as it should be.
 Then the set
 (of cardinality ${n(\sigma)-1}\choose {k-1}$ ) 
 of all $\zeta_{\sigma(a)}$
is a basis
of the subspace of $\tilde H_p(0,U)$ generated by $\zeta_{\sigma}$ with
$\rk(\sigma)=1$.
\end{remark}

{\bf Example.}

\medskip
Let us consider the case where $n=6$ and $k=3$.  
The results for the additive structure of $H^*(M_{6,3})$
can be expressed as the following table (cf. \cite{BW}, p. 311).

$$\matrix{
 \dim U& |\sigma|& p&
 {\rm rank}(\sigma)&
                 s(\sigma)& n(\sigma)& q
& {\rm dim}\tilde H^q(M)\cr
&&&&&&&\cr
 2& 1& -1& 1& 1
       & 3& 3& 20\cr
 3& 2& 0& 1& 1
       & 4& 4& 45\cr
 4& 3& 1& 1& 1
       & 5& 5& 36\cr
 5& 4& 2& 1& 1
       & 6& 6& 10\cr
 4& 2& 0& 2& 2
       & 6& 6& 10\cr
 5& 3& 1& 2& 1
       & 6& 7& 10\cr
}$$

In this table, $U=\bigvee(\sigma)$, $p=|\sigma|-2$ is the dimension of
the homology of $[0,U]$ generated by $\zeta_{\sigma}$ and $q$ is the respective
dimension of the cohomology $H^*(M)$ with $M=M_{6,3}$. Notice that $H^6(M)$ 
is represented as $H^6(M)=H^6_1(M)\oplus H^6_2(M)$
where $H^6_1(M)$ is generated by $\zeta_{\sigma}$ with $|\sigma|=4$ and
$\rk(\sigma)=1$ and $H^6_2(M)$ is generated by $\zeta_{\sigma}$ with
$|\sigma|=2$ and $\rk(\sigma)=2$.

As a ring, $H^*(M)$ is generated by 
$H^i(M),\ i=3,4,5,$ and $H^6_1(M)$. Using Remark \ref{rank1} it is easy to
exhibit a concrete set of algebar generators $\zeta_{\sigma}$.
The multiplication $\nu$ is not zero only on
$H^{3,3}=H^3(M)\otimes H^3(M)$ and $H^{3,4}=H^3(M)\otimes H^4(M)$. 
Due to Theorem \ref{k-equal}, $\nu(\zeta_{\sigma}\otimes\zeta_{\tau})\not=0$
on $H^{3,3}$ only if $b_A\cap b_B=\emptyset$ where $\sigma=\{A\},\
\tau=\{B\}$ and $b_A,b_B$ are respective non-trivial blocks of the atoms $A$
and $B$. This  and the skew commutativity of $\nu$ give 10 linearly
independent products $\zeta_{\sigma\cup\tau}$
in $H^6_2(M)$ that generate this whole space.

On $H^{3,4}$, consider the product $\pi=\zeta_{\sigma\cup\tau}$ where
$\sigma=\{A\}$ and $\tau=\{B_1,B_2\}$ with $|b_{B_1}\cap
b_{B_2}|=2$. Then Theorem \ref{k-equal} implies that
$\pi\not=0$ only if either $b_{B_1}$ or
$b_{B_2}$ is disjoint with $b_A$ and in this case $\pi$ does not depend on the
other atom in $\tau$. Again it is clear that there are 10 linearly
independent products and they generate the whole space $H^7(M)$.

It may be instructive to exhibit a independent set of atoms $\sigma$ such that
$\zeta_{\sigma}=0$. The simplest examples are given by $\sigma=\{B_1,B_2\}$
with $|b_{B_1}\cap b_{B_2}|=1$. For instance, if
$\sigma=\{\{1,2,3\},\{1,4,5\}\}$ 
where the atoms are given by their non-trivial blocks then $\omega=\sigma\cup
\{\{1,2,4\}\}$ is a relater of $[0,\bigvee(\sigma)]$ and the relation
$\partial\omega$ gives $\zeta_{\sigma}=0$.
For a less obvious example one can take 
$\sigma=\{\{1,2,3\},\{3,4,5\},\{3,5,6
\}\}$. Then 
$\rk(\sigma)=1$ and $3=|\sigma|\not=n(\sigma)-\rk(\sigma)(k-2)-s(\sigma)=4$.
On the other hand, if $\omega=\sigma\cup\{\{1,3,5\}\}$ then it is a relater of
$\Pi_{6,3}$ and the relation
$r_{\omega}$ gives $\zeta_{\sigma}=0$.


\bigskip
\begin{thebibliography}{99}

\bibitem{Ar} V. I. Arnold, The cohomology ring of the colored braid group, Mat.
Zametki {\bf 5} (1969), 227-231 (Math. Notes {\bf 5} (1969), 138-140).


\bibitem{BWa} A. Bj\"orner and J. Walker, A homotopy complementation formula
for partially ordered sets, European J. Combin {\bf 4} (1983), 11-19.


\bibitem{BW} A. Bj\"orner and V. Welker, The homology of ``$k$-equal''
manifolds and related partition lattices, Advances in Math. {\bf 110}
(1995), 277-313.


\bibitem{Br} E. Brieskorn, Sur les groupes de tresses, in S\'eminare Bourbaki
1971/72, Lecture Notes in Math {\bf 317}, Springer Verlag, 1973, pp. 21-44.


\bibitem{DCP} C. De Concini and C. Procesi, Wonderful models of subspace
arrangements, Selecta Mathematica, New Series, 1(1995), 459-494.


\bibitem{Fe} E. M. Feichtner, Cohomology algebras of subspace arrangements and of
classical configuration spaces, Cuvillier Verlag G\"ottingen, 1997 (Doctors
Dissertation at TU, Berlin). 


\bibitem{Fo} J. Folkman, The homology groups of a lattice, J. Math. and Mech.
{\bf 15} (1966), 631-636.


\bibitem{Ga} G.Gaiffi, Blow-ups and cohomology bases for De Concini -Procesi
models of subspace arrangements, preprint, 1996.


\bibitem{GMP} M.Goresky and R.Mac Pherson, ``Stratified Morse Theory'',
Part III, Springer
Verlag, 1988.


\bibitem{Mo} J.Morgan, The algebraic topology of smooth algebraic varieties,
Publ. Math. IHES {\bf 48} (1978), 137-204.


\bibitem{Mu} J.Munkres, ``Elements of algebraic topology'', Addison-Wesley,
Menlo Park, CA, 1984.


\bibitem{OS} P.Orlik and L.Solomon, Combinatorics and Topology of Complements of
Hyperplanes, Invent. Math. {\bf 56} (1980), 167-189.


\bibitem{OT} P.Orlik and H.Terao, ``Arrangements of hyperplanes'',
Springer-Verlag, Berlin, 1992.


\bibitem{Qu} D. Quillen, Homotopy properties of the poset of nontrivial
$p$-subgroups of a group, Advances in Math. {\bf 28} (1978), 101-128.


\bibitem{Yu} S.Yuzvinsky,
Cohomology bases for the De Concini - Procesi models of hyperplane
arrangements and sums over trees, Invent. Math. {\bf 127} (1997), 319-335. 


\bibitem{ZZ} G.Ziegler and R.\v Zivaljevi\'c. Homotopy types of subspace
arrangements via diagrams of spaces, Math. Ann. {\bf 295} (1993), 527-548.

 \end{thebibliography}
\end{document}